\UseRawInputEncoding
\documentclass[11pt,a4paper]{article}
\usepackage[numbers,sort&compress]{natbib}
\usepackage[colorlinks,
            linkcolor=blue,
            anchorcolor=green,
            citecolor=magenta
            ]{hyperref}
\usepackage{mathrsfs,amssymb,amsfonts}
\usepackage{graphicx}
\usepackage{amsmath,amssymb,latexsym,color}
\usepackage[mathscr]{eucal}
\usepackage[numbers]{natbib}
\usepackage{bm}
\usepackage[top=1in, bottom=1in, left=1.25in, right=1.25in]{geometry}
\usepackage{indentfirst}
\usepackage{caption}
\usepackage{subfigure}
\makeatletter

\newcommand{\Rmnum}[1]{\expandafter\@slowromancap\romannumeral #1@}
\makeatother
\newtheorem{theorem}{Theorem}[section]

\newtheorem{lemma}[theorem]{Lemma}
\newtheorem{remark}[theorem]{Remark}

\newtheorem{corollary}[theorem]{Corollary}
\newtheorem{example}{Example}

\newcommand{\qed}{\hfill\Box\medskip}
\usepackage{geometry}
\begin{document}
%\begin{CJK*}{GBK}{song}
 \setlength{\baselineskip}{14pt}
\renewcommand{\abovewithdelims}[2]{
\genfrac{[}{]}{0pt}{}{#1}{#2}}
%%%%%%%%%%%%%%%%%%%%%%%%%%%%%%%%%%%%%%%%%%%%%%%%%%%%%%%%%%%%%%%%%%%%%%%%%%%%%%%%%%%%%%%%
%%%%%%%%%%%%%%%%%%%%%%%%%%%%%%%%%%%%%%%%%%%%%%%%%%%%%%%%%%%%%%%%%%%%%%%%%%%%%%%%%%%%%%%%

\title{\bf Component Edge Connectivity of Hypercube-like Networks \footnote{This research is supported partially by the General Project of Hunan Provincial Education Department of China (19C1742), Hu Xiang Gao Ceng Ci Ren Cai Ju Jiao Gong Cheng-Chuang Xin Ren Cai (No. 2019RS1057), the Project of Scientific Research Fund of Hunan Provincial Science and Technology Department (No. 2018WK4006).}
%\footnote{For the title, try not to use more than 3
%lines. Typeset the title in 10 pt Roman, boldface with the first letter of
%important words capitalized.}
}

%\title{\bf Edge-fault-tolerant strong Menger edge connectivity of hypercubes and folded hypercubes\\
%\footnotesize{ This research is supported by the National Natural Science Foundation of China (11571044, 61373021),  the Fundamental Research Funds for the Central University of China.} }

\author{Dong Liu \quad   Pingshan Li
 \footnote{Corresponding author. \newline {\em E-mail address:} lips@xtu.edu.cn (Pingshan Li).} \quad Bicheng Zhang \\
 {\footnotesize \em Key Laboratory of Intelligent computing$\And$Information processing of Ministry of Education}\\
{\footnotesize   \em
School of  Mathematical and Computational Sicence, Xiangtan University, Xiangtan, Hunan 41105, PR China }}
 \date{}
 \date{}
 \maketitle

\begin{abstract}
As a generalization of the traditional connectivity, the  $g$-component edge connectivity $c\lambda_g(G) $ of a non-complete graph $G$ is the minimum number of edges to be deleted from the graph $G$ such that the resulting graph has at least $g$  components. Hypercube-like networks (HL-networks for short) are obtained by manipulating some pairs of edges in hypercubes, which contain several famous interconnection networks such as twisted cubes, M$\ddot{o}$bius cubes,  crossed cubes, locally twisted cubes. In this paper, we determine the $(g+1)$-component edge connectivity of the $n$-dimensional HL-networks for $g\le 2^{\left\lceil \frac{n}{2} \right\rceil}$, $n \ge 8$.


\medskip
%\noindent {\em 2010 MSC:} 05C25; 05C15.
\end{abstract}

\section{Introduction}
As we all know, the interconnection network is an important part of the multiprocessor system. For convenience, we usually model the interconnection network by a graph, with vertices representing processors and edges representing links between processors. There are many parameters to evaluate the reliability of a network. The traditional connectivity $\kappa(G)$ of the graph $G$ is the most crucial one among them.  Generally speaking, the larger the $\kappa(G)$ of the network is, the better its reliability is. However, the traditional connectivity only indicates when the network will break but does not further depict the properties of  components, which makes it impossible to accurately evaluate the reliability of the network. In order to further describe the properties of  components, a lot of more precise connectivity  have been proposed,  such as extra connectivity \rm{\cite{2010Conditional}}, super connectivity  \rm{\cite{1981Connectivity}} and restricted connectivity \rm{\cite{1989Generalized}}. In 1984, Chartrand et al. \rm{\cite{1984Generalized}} and Sampathkumar \rm{\cite{1984Connectivity}} respectively introduced $g$-component connectivity $c\kappa_g(G)$ and $g$-component edge connectivity $c\lambda_g(G)$ of the graph $G$.

A $g$-component (edge) cut  of a non-complete $G$ is a vertex (edge) set to be deleted from the graph $G$ such that the resulting graph has at least $g$  components. The $g$-component (edge) connectivity of $G$, written $c\kappa_g(G) $ ($c\lambda_g(G) $), is the minimum size of the  $g$-component (edge) cut of $G$. Obviously, the component (edge) connectivity is a generalization of the traditional connectivity and $\kappa(G)=c\kappa_2(G)$ ($\lambda(G)=c\lambda_2(G)$). More  importantly, compared with the traditional connectivity, the component (edge) connectivity can be better satisfied in practical applications. Thus, many researches have been focused on the component connectivity of some famous networks in recent years. (see, for example {\rm \cite{2015Connectivity,2017Structural,2019Characterizations,2019Conditional}}). But there are relatively few papers about component edge connectivity. Zhao et al. \rm{\cite{Zhao2018Component}} determined the $(g+1)$-component edge connectivity of hypercubes for $g\le 2^{\left[\frac{n}{2} \right]}$, $n \ge 7$. Based on this nice work, we shall be centered on the $(g+1)$-component edge connectivity of hypercube-like networks for $g\le 2^{\left\lceil \frac{n}{2} \right\rceil}$, $n \ge 8$.

As a popular topology for the design of the multiprocessor system, hypercubes have many excellent properties including  symmetry, relatively diameter, good connectivity and recursive scalability. Improving these properties has led to the evolution of hypercube variants.  A lot of  hypercube variants,  such as twisted cubes \rm{\cite{2002Fault}},  crossed cubes \rm{\cite{Efe1991A}} and M$\ddot{o}$bius cubes  \rm{\cite{1995The}}, have been introduced successively. Many of the properties of these networks are identical with that of the hypercubes. In particular, their diameter is shorter than that of the hypercubes. In order to conduct a unified study on these variants, Vaidya et al. proposed hypercube-like networks (in brief, HL-networks), which contain all of the networks mentioned above. Hence, HL-networks  attracted considerable attention in the past (see, for example {\rm \cite{2013,Jin2017On}}). This class of networks are sometimes called BC-networks \rm{\cite{2003BC}}. we will use the term HL-networks in this paper.

\begin{figure}[htbp]
\centering	
	\subfigure[]
	{
		\begin{minipage}{6cm}
			\centering
			\includegraphics[scale=0.4]{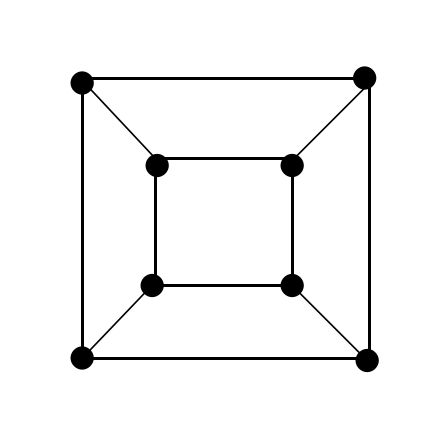}
		\end{minipage}
	}
	\subfigure[]
	{
		\begin{minipage}{6cm}
			\centering
			\includegraphics[scale=0.4]{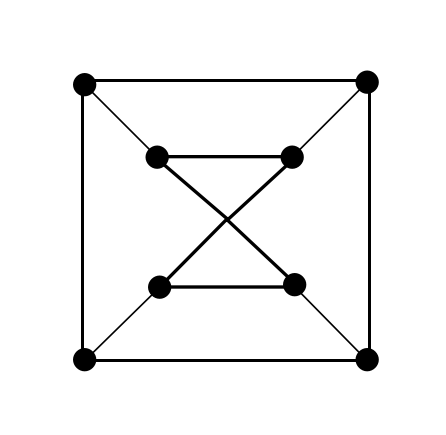}
		\end{minipage}
	}
	\caption{Two 3-dimensional HL-networks: (a) $Q_3$\ \  (b) $G(8,4)$}
    \label{f1}
\end{figure}

The recursive definition of the HL-networks is as follows:

\ \ $HL_1=\{K_1\}$ and $HL_{n}=\{G_{n-1}\oplus G^*_{n-1}|G_{n-1},G^*_{n-1}\in HL_{n-1}\}$,\\
where the symbol $``\oplus"$ represents the perfect matching operation that connects $G_{n-1}$ and $G^*_{n-1}$ using some disjoint edges.
It's easy to get that $HL_0=\left\{ K_1\right\}$,  $HL_1=\left\{ K_2\right\}$, $HL_2=\left\{ C_4\right\}$, and $HL_3=\left\{ Q_3,\  G(8,4)\right\}$, where $C_4$ is a cycle of length 4, and $Q_3$ and $G(8,4)$ are shown in Figure \ref{f1}. The $n$-dimensional HL-network $G_n$ is $n$ regular, and it has $2^n$ vertices and $n2^{n-1}$ edges.

Next, we will introduce some notations and definitions which will be used in this paper. Let $G=(V(G),E(G))$ be a graph. The size of G is the number of edges of $G$.  The degree of $v$, denoted by $d_{G}(v)$, is the number of edges incident to $v$ in $G$. For a subset $X$ of $V(G)$,  $G[X]$ is the subgraph induced by  $X$. If $D_i$ and $D_j$ are two disjoint subgraphs of $G$, then we use $E(D_i,D_j)$ to denote the set of edges between the subgraphs $D_i$ and $D_j$. Similarly, we use $E(v, D_i)$ to denote the set of edges between  $v$ and $D_i$, where $v\in V(G)$ but $v \notin V(D_i)$. For any $v\in V(G)$,  %we define $G-v=G[V(G)/v]$ and $G-D_i=G[V(G)/V(D_i)]$.
$G-v$ denotes a subgraph obtained by deleting  $v$ and edges incident to $v$. Similarly, $G-D_i$ denotes a subgraph obtained by removing all vertices in $V(D_i)$ and all edges incident to vertices of $V(D_i)$, where $D_i$ is a subgraph of $G$.

The rest of this paper is organized as follows. In section $2$, we shall determine the maximum size of the subgraph induced by $g$ vertices in  HL-networks. In section $3$, we shall determine
the $(g+1)$-component edge connectivity of HL-networks by applying the results of section 2.

\section{The maximum size of the subgraph induced by $g$ vertices  }
It is a classical problem to determine the  size (the number of edges) of the subgraph that satisfies some given property in a graph. For instance, Erd$\ddot{o}$s has studied an interesting problem in \rm{\cite{1990A}}: What is the maximum size of the subgraph without cycles of length 4 in hypercubes. In this section, we shall determine the  maximum size of the subgraph induced by $g$ vertices in HL-networks. Furthermore, its application in the next section helps us find the $(g+1)$-component edge connectivity of HL-networks.

Let $e_g$ be the maximum size of the subgraph induced by $g$ vertices in the $n$-dimensional HL-network $G_n$, that is, $e_g=$ $\max\left\{|E( G_n[X])|: X \subseteq V(G_n)\ {\rm and}\ |X|=g\right\}$. By the mathematical principle  of converting decimal digit to binary digit, any integer $g$ can be written as the sum of the exponents of 2, that is, $g=\begin{matrix} \sum_{i=0}^s 2^{t_i}\end{matrix} $, where $t_0 = \left\lfloor \log_2 g\right \rfloor$, $t_i = \left\lfloor \log_2(g- \begin{matrix} \sum_{r=0}^{i-1} 2^{t_r}\end{matrix})\right\rfloor $ for $i \ge 1$.  Li and Yang have determined $e_g$ of hypercubes in \rm{\cite{2013Bounding}}. In order to avoid confusion, we use $e_g(Q_n)$ to denote $e_g$ of hypercubes, where $e_g(Q_n)=\begin{matrix} \sum_{i=0}^s t_i2^{t_i-1}\end{matrix}+\begin{matrix} \sum_{i=0}^s i2^{t_i}\end{matrix}$. Furthermore, the function $e_g(Q_n)$ has the following property.

\begin{lemma}\rm{\cite{2013Bounding}}\rm \cite{Zhao2018Component}\label{le2.1}
If $g_0\le g_1$, then $e_{g_0+g_1}(Q_n) \ge e_{g_0}(Q_n)+e_{g_1}(Q_n)+g_0$.
\end{lemma}
\begin{figure}[t]
    \centering
    \includegraphics[scale=0.4]{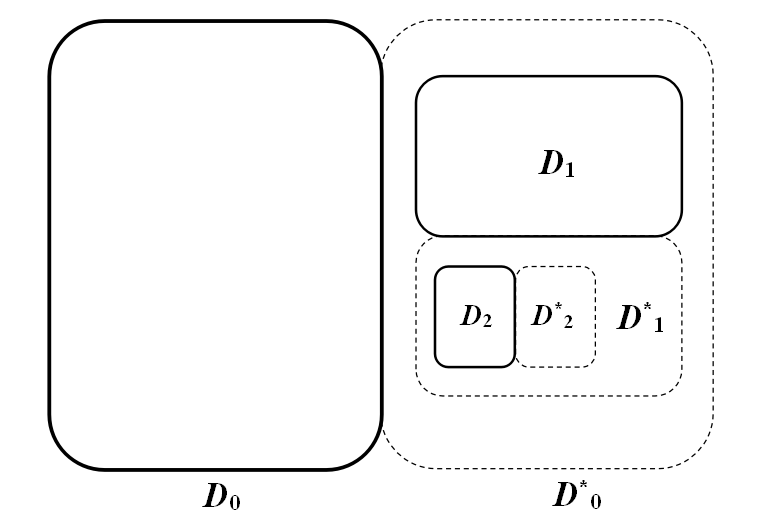}
    \caption{The HL-networks $D_i$ and $D^*_i$ when {\em s}=2}
    \label{f2}
    \end{figure}

Next, we shall give an algorithm to find the subgraph $H$ induced by $g$ vertices in the $n$-dimensional HL-network $G_n$ such that $|E(H)|= \begin{matrix} \sum_{i=0}^s t_i2^{t_i-1}\end{matrix}+\begin{matrix} \sum_{i=0}^s i2^{t_i}\end{matrix}$. \\
\rule[0.25\baselineskip]{\textwidth}{1pt}\\
Algorithm-$MS$

{\bf Input}: An $n$-dimensional network $G_n$, two integer $g$ (suppose $g=\begin{matrix} \sum_{i=0}^s 2^{t_i}\end{matrix}$) and $s$.

{\bf Output}: A vertex set $V=V(D_0)\cup V(D_1)...\cup V(D_s)$ and a subgraph $H=G_n[V]$.

Initialization: $i\longleftarrow 0$, $V \longleftarrow  \emptyset$, where $\emptyset$ is an empty set.

Iteration: As long as $i\leq s$, we take a $t_i+1$-dimensional subcube $T_i$ from $D^*_{i-1}$ $(D^*_{-1}=G_n)$. Then $T_{i}$ can be written as $D_i\oplus D^*_i$, where $D_i$ and $D^*_i$ are two $t_i$-dimensional subcubes. Add $V(D_i)$ to $V$, $i \longleftarrow i+1$.\\
\rule[0.25\baselineskip]{\textwidth}{1pt}

For ease of understanding, we list $D_i$ as follows (see Figure \ref{f2} for $s=2$):

\begin{equation*}\begin{split}
D_0 \ \ \ & (t_0\text{-dimensional subcube}) ; \\
D_1 \ \ \ & (t_1\text{-dimensional  subcube taken from}\   D^*_0);\\
D_2 \ \ \ & (t_2\text{-dimensional subcube taken  from}\ D^*_1);\\
......\\
D_s \ \ \ & (t_s\text{-dimensional subcube taken from}\  D^*_{s-1}).
\end{split}
\end{equation*}

We need to verify that the subgraph $H$ found by the algorithm-$MS$  that satisfies $|V(H)|=g$ and $|E(H)|=\begin{matrix} \sum_{i=0}^s t_i2^{t_i-1}\end{matrix}+\begin{matrix} \sum_{i=0}^s i2^{t_i}\end{matrix}$.

Note that $H=G_n[V(D_0)\cup V(D_1)...\cup V(D_s)]$ and $|V(D_i)|=2^{t_i}$. Thus,
\begin{center}
 $V(H)=\sum_{i=0}^s|V(D_i)|=\sum_{i=0}^s2^{t_i}=g$.
\end{center}

By the choice of $D_i$, one has $|E(D_i,D_j)|=|V(D_j)|=2^{t_j}$ for $i<j$. Thus,
\begin{center}
$E(H)=\sum_{i=0}^s |E(D_i)|+ \sum_{i=0}^{s-1}\sum_{j=i+1}^s|E(D_i\ ,D_j)|=\sum_{i=0}^s t_i2^{t_i-1}+ \sum_{i=0}^s i2^{t_i}$.
\end{center}
\begin{theorem}\label{TH2.2}
For any $G_n\in HL_n$, suppose that $ X \subseteq V(G_n)$ with $|X|=g=\begin{matrix} \sum_{i=0}^s 2^{t_i}\end{matrix}$, then $e_g=\begin{matrix} \sum_{i=0}^s t_i2^{t_i-1}\end{matrix}+\begin{matrix} \sum_{i=0}^s i2^{t_i}\end{matrix}$.
\end{theorem}

\noindent{\bf Proof: }For convenience, we define $f(g)=\begin{matrix} \sum_{i=0}^s t_i2^{t_i-1}\end{matrix}+\begin{matrix} \sum_{i=0}^s i2^{t_i}\end{matrix}$. Firstly, we shall  prove $|E(G_n[X])|\le f(g)$ for any $X\in V(G_n)$ with $|X| =g=\begin{matrix} \sum_{i=0}^s 2^{t_i}\end{matrix}$ by induction on $n$. Clearly, the result holds for $n=2$ since $G_2 \cong C_4$. So we assume $n \ge 3$. Note that $G_n=G_{n-1}\oplus G^*_{n-1}$. Let $X_0=X \cap G_{n-1}$,  $X_1=X \cap G^*_{n-1}$, $|X_i|=g_i$ for $i=0, 1$ and $g=g_0+g_1$. Without loss of generality, we suppose $g_1\geq g_0$. Note $f(g)=e_g(Q_n)$. By induction hypothesis, we have that
\begin{equation*}\begin{split}\label{eq2}
 |E(G_n[X])| &\le |E(G_n[X_0])|+|E(G_n[X_1])|+g_0\\
 &=|E(G_{n-1}[X_0])|+|E(G^*_{n-1}[X_1])|+g_0\\
 &\le f(g_0)+f(g_1)+g_0\\
 &=e_{g_0}(Q_n)+e_{g_1}(Q_n)+g_0\\
 &\le e_{g_0+g_1}(Q_n)=f(g_0+g_1)=f(g).
\end{split}
\end{equation*}
The third inequality holds because of Lemma \ref{le2.1}.

Using the algorithm-$MS$, we can  find a subgraph $H$ such that $|V(H)|=g$ and $|E(H)|=\begin{matrix} \sum_{i=0}^s t_i2^{t_i-1}\end{matrix}+\begin{matrix} \sum_{i=0}^s i2^{t_i}\end{matrix}=f(g)$. Thus, $e_g=f(g)=\begin{matrix} \sum_{i=0}^s t_i2^{t_i-1}\end{matrix}+\begin{matrix} \sum_{i=0}^s i2^{t_i}\end{matrix}$.
$\qed$

Here are some properties of the function $e_g$.
\begin{lemma}\label{le2.3}
If $g\le 2^{n-2}$, then $(n-2)g-2e_g\ge 0$.
\end{lemma}
\noindent{\bf Proof: }If $g\le 2^{n-2}$, then we can take a subgraph $H$ from an  ($n-2$)-dimensional HL-network $G_{n-2}$ such that $|V(H)|=g$ and $|E(H)|=e_g$. Note that $G_{n-2}$ is an ($n-2$)-regular graph. Thus, $(n-2)g-2e_g\geq 0$.
$\qed$
\begin{lemma}\label{le2.4}
$e_{i+1}=e_i+s+1$, where $i=2^{t_0}+2^{t_1}+...+2^{t_s}$.
\end{lemma}
\noindent{\bf Proof:}
If $t_s>0$, then $i+1=2^{t_0}+2^{t_1}+...+2^{t_s}+2^0 $. Note that $e_i=\begin{matrix} \sum_{i=0}^s t_i2^{t_i-1}\end{matrix}+\begin{matrix} \sum_{i=0}^s i2^{t_i}\end{matrix}$. We have that
\begin{center}
$e_{i+1}=\begin{matrix} \sum_{i=0}^s t_i2^{t_i-1}\end{matrix}+0\cdot2^{0-1}+\begin{matrix} \sum_{i=0}^s i2^{t_i}\end{matrix}+(s+1)\cdot2^0=e_i+s+1$.
\end{center}

Otherwise, $t_s=0$. Then there is an integer $r$ ($0\le r \le s $) such that $t_j=s-j$ for all $r \leq j\leq s$. In other words, $i=2^{t_0}+...+2^{t_{r-1}}+2^{s-r}+2^{s-r-1}+...+2^1+2^0$. Then $i+1=2^{t_0}+...+2^{t_{r-1}}+2^{s-r}+2^{s-r-1}+...+2^1+2^0+2^0=2^{t_0}+...+2^{t_{r-1}}+2^{s-r+1}$.

Using the algorithm-$MS$, we take two subgraph $H_1$ and $H_2$ from $G_n$, where $H_1$ and $H_2$ are induced by $i$ and $i+1$ vertices respectively. Clearly, $|E(H_1)|=e_i$ and $|E(H_2)|=e_{i+1}$. we can assume that
\begin{equation}\begin{split}\label{eq1}
& H_1=G_n\left[V(D_0)\cup ...\cup V(D_{r-1})\cup V(D_r)\cup V(D_{r+1})\cup ...\cup V(D_s)\right];\\
& H_2=G_n\left[V(D_0)\cup ...\cup V(D_{r-1})\cup V(D'_r)\right].
\end{split}
\end{equation}
where $D_i$ is a $t_i$-dimensional subcube, $D'_r$ is an $(s-r+1)$-dimensional subcube.

In fact, $|V(H_1)|=i$ and $|V(H_2)|=i+1$. In other words, $H_1$ has one less vertex than $H_2$. Let the vertex be $v$. From (\ref{eq1}),  we can see that the first $r$ subcubes $D_i$ in $H_1$ and $H_2$ are the same. Thus, one has $v\subseteq V(D'_r)$.  Clearly,  $d_{D'_r}(v)=s-r+1$. In addition, there is only one edge between $v$ and every $D_i$ in $H_2$ for $0\leq i\leq r-1$. Thus, $d_{H_2}(v)=r+(s-r+1)=s+1$. We have that $|E(H_2)|=|E(H_1)|+d_{H_2}(v)$, that is, $e_{i+1}=e_i+s+1$.
$\qed$
\begin{lemma}\label{le2.5}
If $i\le j$, then $e_{i+1}+e_j\le e_{i+j}$.
\end{lemma}
\noindent{\bf Proof: } Suppose $i=2^{t_0}+...+2^{t_s}$, then $e_{i+1}=e_i+s+1$ by Lemma \ref{le2.4}. Note that $e_g(Q_n)=e_g$. By the Lemma \ref{le2.1}, $e_{i+j}\ge e_i+e_j+i=e_{i+1}+e_{j}+i-(s+1)\ge e_{i+1}+e_{j}$.
$\qed$

\section{Component edge connectivity of HL-networks}

In this section, we shall apply Theorem \ref{TH2.2} to determine the $(g+1)$-component edge connectivity of HL-networks.

  For any $X\subseteq V(G_n)$, we use $E_X$ to denote a set of edges in which each edge has exactly one endpoint in $X$.
\begin{lemma}\label{le3.1}
For any $G_n\in HL_n$,  let $X\subseteq V(G_n)$ with $|X|=g=\begin{matrix} \sum_{i=0}^s 2^{t_i}\end{matrix}$. Then $|E_X|\ge ng-2e_g$. Moreover,  $ng-e_g$ is strictly increasing (respect to $g$) for $g\le 2^{\left \lceil  \frac{n}{2} \right \rceil  },\ n\ge 2$.
\end{lemma}
\noindent{\bf Proof: }Note that $G_n$ is an $n$-regular graph. We find that $|E_X|=n|X|-2|E(G_n[X])|\ge ng-2e_g$. According to Lemma \ref{le2.4}, if $g\le 2^{\left \lceil  \frac{n}{2} \right \rceil  }$ and $n\ge 2$, then $e_{g+1}-e_g=s+1<n$. Hence,  we have  that $(n(g+1)-e_{g+1})-(ng-e_g)=n-(e_{g+1}-e_g)>0$, which indicate that $ng-e_g$ is strictly increasing for $g\le 2^{\left \lceil  \frac{n}{2} \right \rceil },\ n\ge 2$.
$\qed$
\begin{theorem}
For any $G_n\in HL_n$, $c\lambda_{g+1}(G_n)=ng-e_g $ for $g \le 2^{\left \lceil \frac{n}{2} \right \rceil}, \ n\ge 8$.\label{th3.2}
\end{theorem}
\noindent{\bf Proof: }We first prove that $c\lambda_{g+1}(G_n)\le ng-e_g $ by constructing a $(g+1)$-component edge-cut $F$ such that $|F|=ng-e_g$.  Using the algorithm-$MS$, we can get a subgraph $H$ such that $|V(H)|=g$ and $|E(H)|=e_g$. Let $F=E_{V(H)}\cup E(H)$, then we have that $|F|=ng-e_g$. Moreover, $G_n-F $ has at least $g+1$ components. Thus, $c\lambda_{g+1}(G_n)\le ng-e_g$.
%\begin{figure}[htbp]
%\centering	
	%\subfigure[]
	%{
		%\begin{minipage}{8cm}
			%\centering
			%\includegraphics[scale=0.4]{F1.png}
            %\label{fa}
		%\end{minipage}
	%}
	%\subfigure[]
	%{
		%\begin{minipage}{8cm}
			%\centering
			%\includegraphics[scale=0.4]{F2.png}
            %\label{fb}
		%\end{minipage}
	%}
	%\caption{The edges between the components}
%\end{figure}

Next, we shall prove that $c\lambda_{g+1}(G_n)\ge ng-e_g $.  Assume that $F$ is  the smallest $(g+1)$-component edge-cut and $G_n-F$ has exactly $g+1$ components. Denote the $g+1$ components in $G_n-F$ by $C_1,C_2,...,C_{g+1}$. Without loss of generality, we can assume that $|V(C_1)|\le |V(C_2)|\le...\le |V(C_{g+1})|$.\\
\noindent{\bf Case 1: }Suppose that $|V(C_{g+1})|<2^{n-2}$.

In this case, $|V(C_i)|<2^{n-2}$ for all $i$. Note that $\sum^{g+1}_{i=1}|V(C_i)|=2^n$. Then we can find an integer $j$ such that $ \sum_{i=1}^{j}|V(C_i)|<2^{n-2}$ but $2^{n-2}\le\sum_{i=1}^{j+1}|V(C_i)|< 2^{n-1}$. Let $X=\bigcup_{i=1}^{j+1}V(C_i)$ with $|X|=m=\begin{matrix} \sum_{i=0}^s 2^{t_i}\end{matrix}$. Clearly, $2^{n-2}\le m=\begin{matrix} \sum_{i=0}^s 2^{t_i}\end{matrix} < 2^{n-1}$. It follow that $t_0=n-2$. We suppose that $m'=m-2^{t_0}$, then $m'<2^{n-2}$. By Lemma \ref{le2.3}, we have that $(n-2)m'-2e_{m'}\ge 0$. Combining with Lemma \ref{le3.1}, we have that
\begin{center}$\begin{array}{rl}
|E_X|\ge & nm-2e_m \vspace{1.0ex}\\
= & n2^{t_0}+n(2^{t_1}+...+2^{t_s})-t_02^{t_0}-(\begin{matrix} \sum_{i=1}^s t_i2^{t_i}\end{matrix}+\begin{matrix} \sum_{i=1}^s 2\cdot i\cdot2^{t_i}\end{matrix})\vspace{1.0ex} \\
 = & (n-t_0)2^{t_0}+(n\begin{matrix} \sum_{i=1}^s 2^{t_i}\end{matrix}-(\begin{matrix} \sum_{i=1}^s t_i2^{t_i}\end{matrix}+\begin{matrix} \sum_{i=1}^s 2\cdot i\cdot2^{t_i}\end{matrix})) \vspace{1.0ex} \\
= & (n-t_0)2^{t_0}+nm'-2e_{m'}-2m' \vspace{1.0ex} \\
\ge & (n-t_0)2^{t_0}\vspace{1.0ex} \\
= & (n-t_0)\cdot 2^{(t_0-\left \lceil \frac{n}{2} \right \rceil)}\cdot 2^{\left \lceil\frac{n}{2} \right \rceil}. \\
\end{array}$
\end{center}
 Note that $t_0=n-2$.  Then $(n-t_0)\cdot 2^{(t_0-\left \lceil \frac{n}{2} \right \rceil)}>  \frac{3n}{4}$ for $n\ge 8$.  we have that
\begin{center}
$|F|\ge |E_X|>\frac{3n}{4}  \cdot 2^{\left \lceil \frac{n}{2} \right \rceil}\ge n\cdot 2^{\left \lceil \frac{n}{2} \right \rceil}-e_{2^{\left \lceil \frac{n}{2} \right \rceil}}\ge ng-e_g$.
\end{center}
The last inequality holds, because $ng-e_g$ is strictly increasing (respect to $g$) for $g\le 2^{\left \lceil \frac{n}{2} \right \rceil }$ by Lemma \ref{le3.1}.

\noindent{\bf Case 2: }Suppose that $|V(C_{g+1}|)\ge 2^{n-2}$.

Let $|V(C_i)|=m_i$ for $1\le i\le g$ and  let $m=\begin{matrix} \sum_{i=1}^g m_i \end{matrix}$.  Set $X$=$\bigcup_{i=1}^{g}V(C_i)$.
If $m\ge 2^{n-2}$, we can  get that $ |F|\ge |E_X|> ng-e_g$ by using the proof similar to case 1 . So we  assume  $m < 2^{n-2}$. If $m=g$, then the component $C_i$ ($1\le i\le g $) is an isolated vertex. Thus, $ |F|\ge ng-e_g$. If $m>g$, then let $F'=F\cap E(G_n[X])$. We have  that $ |F|\ge |\bigcup_{i=1}^{g}E_{V(C_i)}|=|E_{V(C_1)}|+...+|E_{V(C_g)}|-|F'|$. Since $ |E_{V(C_i)}|\geq nm_i-2|E(C_i)|$ and $|F'|\le e_m-|\bigcup_{i=1}^{g}E(C_i)|$, we have that
\begin{center}$\begin{array}{rl}
 |F|\ge & |\bigcup_{i=1}^{g}E_{V(C_i)}|=|E_{V(C_1)}|+...+|E_{V(C_g)}|-|F'| \vspace{1.0ex} \\
\ge & \begin{matrix} \sum_{i=1}^g (nm_i- 2|E(C_i)|) \end{matrix}-e_m+\begin{matrix} \sum_{i=1}^g |E(C_i)| \end{matrix} \vspace{1.0ex} \\
= & nm-e_m- \begin{matrix} \sum_{i=1}^g |E(C_i)| \end{matrix}\vspace{1.0ex} \\
\ge & nm-e_m- \begin{matrix} \sum_{i=1}^g e_{m_i} \end{matrix}.\\
\end{array}$
\end{center}

By Lemma \ref{le2.5}, $e_{m_1}+e_{m_2}+...+e_{m_g}\le e_{m_{1}+m_2-1}+e_{m_3}+...+e_{m_g}\le...\le e_{m-g+1}$. Thus, we have that $|F| \ge nm-e_m- \begin{matrix} \sum_{i=1}^g e_{m_i} \end{matrix} \ge nm-e_m- e_{m-g+1}$.

Next, we shall show that $nm-e_m-e_{m-g+1}> ng-e_g$.

\begin{figure}[t]
    \centering
    \includegraphics[scale=0.25]{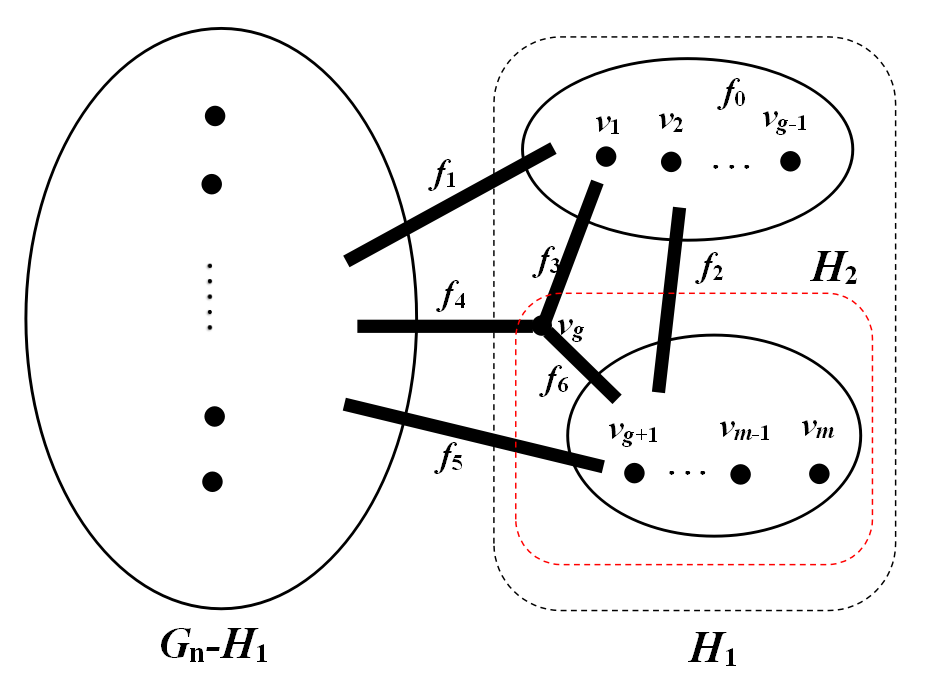}
    \caption{The edges between the components}
    \label{fb}
    \end{figure}

Note that $G_n=G_{n-1}\oplus G^*_{n-1}$. We take a subgraph $H_1$ of $m$ vertices from an ($n-2$)-dimensional subcube in $G_{n-1}$ by using the algorithm-$MS$, since $m<2^{n-2}$. Clearly,  $|E(H_1)|=e_m$.  We take a subgraph $H_2$ of $m-g+1$ vertices from $H_1$ by using the algorithm-$MS$,  since $m-g+1\le m$. So $|E(H_2)|=e_{m-g+1}$.  Use $v_m, v_{m-1},...,v_{g+1}, v_g $ to denote the vertices of $H_2$.  By the proof of Lemma \ref{le2.4}, there exists a vertex in $H_2$, say $v_g$, such that $d_{H_2}(v_g)=s+1$, where $m-g=2^{t_0}+2^{t_1}+...+2^{t_s}$. Use  $v_1, v_2,...,v_{g-1}$ to denote  vertices of $V(H_1)-V(H_2)$. We assume that (see Figure \ref{fb})
\begin{alignat*}{2}
f_0&=|E(H_1-H_2)|; \ \ \ \ \ \ \ \ \ \ \ \ \ \ \ \ \ \ \ \ \ \ \ \  &f_4&=|E(v_g, G_n-H_1)||; \\
f_1&=|E(H_1-H_2,G_n-H_1)|;       &f_5&=|E(G_n-H_1, H_2-v_g)|;\\
f_2&=|E(H_1-H_2,H_2-v_g)|;        &f_6&=|E(v_g, H_2-v_g)|;\\
f_3&=|E(v_g, H_1-H_2)|.
\end{alignat*}
Then we find that
\begin{equation}\label{eq2}
\begin{split}
f_0+f_1+f_2+f_3+f_4+f_5&=(f_1+f_4+f_5)+(f_0+f_2+f_3)\\
&=|E_{V(H_1)}|+(|E(H_1)|-|E(H_2)|)\vspace{1.0ex}\\
&=nm-2|E(H_1)|+(|E(H_1)|-|E(H_2)|)\vspace{1.0ex}\\
&=nm-e_m-e_{m-g+1}.\\
\end{split}
\end{equation}
Let $S=\{v_1, v_2,...,v_g\}$. Note that $|E(G_n[S])|=f_0+f_3\le e_g$, and $|E_S|+2|E(G_n[S])|=(f_1+f_2+f_4+f_6)+2(f_0+f_3)=ng$. Thus,
\begin{equation}\label{eq3}
\begin{split}
f_0+f_1+f_2+f_3+f_4+f_6=(f_1+f_2+f_4+f_6)+(f_0+f_3)\ge ng-e_g.
\end{split}
\end{equation}
Moreover, each vertex of $H_1$ has at least two neighbors out of $H_1$ since $H_1$ be taken from an ($n-2$)-dimensional subcube. Thus, $f_5>|V(H_2-v_g)|= m-g$. Note that $ f_6=d_{H_2}(v_g)=s+1$ and $m-g=2^{t_0}+2^{t_1}+...+2^{t_s}$. It is not difficult to get that $f_5> m-g\ge s+1=f_6$.  Thus,
$nm-e_m-e_{m-g+1}> ng-e_g$ by comparing (\ref{eq2}) with (\ref{eq3}).

To sum up, $c\lambda_{g+1}(G_n)=ng-e_g$.
$\qed$

\begin{corollary}
Let $F$ be a $(g+1)$-component edge-cut of the $n$-dimensional HL-network $G_n$ and $|F|=c\lambda_{g+1}(G_n)$, then $G_n-F$ contains $g$ isolated vertices for $g \le 2^{\left \lceil \frac{n}{2} \right \rceil}, \ n\ge 8$.
\end{corollary}
\noindent{\bf Proof: }
By the proof of Theorem \ref{th3.2}, $|F|=c\lambda_{g+1}(G_n)$ if and only if $m=g$ in case 2, that is, $g$ components are isolated vertices.
$\qed$

\section{Conclusion}
 Component edge connectivity is an generation of the traditional connectivity. In this paper, we studied the $(g+1)$-component edge  connectivity of HL-networks. $c\lambda_{g+1}(G_n)=ng-e_g $ for $g \le 2^{\left \lceil \frac{n}{2} \right \rceil}, \ n\ge 8$. But for $g> 2^{\left \lceil \frac{n}{2} \right \rceil}+1$, The problem has not been solved.

\bibliographystyle{plainnat}      %LaTex Class文件, IEEEtran为给定模板格式定义文件名
\bibliography{reference}                        %ref为.bib文件名

\begin{thebibliography}{18}
\providecommand{\natexlab}[1]{#1}
\providecommand{\url}[1]{\texttt{#1}}
\expandafter\ifx\csname urlstyle\endcsname\relax
  \providecommand{\doi}[1]{doi: #1}\else
  \providecommand{\doi}{doi: \begingroup \urlstyle{rm}\Url}\fi

\bibitem[Baker et~al.(1990)Baker, Bollobas, and Hajnal]{1990A}
A.~Baker, B.~Bollobas, and A.~Hajnal.
\newblock A tribute to paul erd$\ddot{o}$s || some of my favourite unsolved
  problems.
\newblock 10.1017/CBO9780511983917:\penalty0 467--478, 1990.

\bibitem[Bauer et~al.(1981)Bauer, Boesch, Suffel, and
  Tindell]{1981Connectivity}
D.~Bauer, F.~Boesch, C.~Suffel, and R.~Tindell.
\newblock Connectivity extremal problems and the design of reliable
  probabilistic networks.
\newblock \emph{theory \& applications of graphs}, 1981.

\bibitem[Chang and Hsieh(2013)]{2013}
N.~W. Chang and S.~Y. Hsieh.
\newblock {2,3}-extraconnectivities of hypercube-like networks.
\newblock \emph{Journal of Computer and System Sciences}, 79\penalty0
  (5):\penalty0 669¨C688, 2013.

\bibitem[Chartrand et~al.(1984)Chartrand, Kapoor, Lesniak, and
  Lick]{1984Generalized}
G.~Chartrand, S.~F. Kapoor, L.~Lesniak, and D.~R. Lick.
\newblock Generalized connectivity in graphs.
\newblock 1984.

\bibitem[Cheng et~al.(2017)Cheng, Ke, and Shen]{2017Structural}
E.~Cheng, Q.~Ke, and Z.~Shen.
\newblock \emph{Structural Properties of Generalized Exchanged
  HypercubesGeneralized exchanged hypercube}.
\newblock Emergent Computation, 2017.

\bibitem[Cull and Larson(1995)]{1995The}
P.~Cull and S.~M. Larson.
\newblock The mobius cubes.
\newblock \emph{IEEE Trans Computers}, 44\penalty0 (5):\penalty0 647--659,
  1995.

\bibitem[Efe and K.(1991)]{Efe1991A}
Efe and K.
\newblock A variation on the hypercube with lower diameter.
\newblock \emph{Computers IEEE Transactions on}, 40\penalty0 (11):\penalty0
  1312--1316, 1991.

\bibitem[Esfahanian(1989)]{1989Generalized}
A.~H. Esfahanian.
\newblock Generalized measures of fault tolerance with application to n -cube
  networks.
\newblock \emph{Computers IEEE Transactions on}, 38\penalty0 (11):\penalty0
  1586--1591, 1989.

\bibitem[Fan and He(2003)]{2003BC}
J.~X. Fan and L.~Q. He.
\newblock Bc interconnection networks and their properties.
\newblock \emph{Chinese Journal of Computers}, 2003.

\bibitem[Harary(2010)]{2010Conditional}
F.~Harary.
\newblock Conditional connectivity.
\newblock \emph{Networks}, 13, 2010.

\bibitem[Huang et~al.(2002)Huang, Tan, Hung, and Hsu]{2002Fault}
W.~T. Huang, Jjm Tan, C.~N. Hung, and L.~H. Hsu.
\newblock Fault-tolerant hamiltonicity of twisted cubes.
\newblock \emph{Journal of Parallel \& Distributed Computing}, 62\penalty0
  (4):\penalty0 591--604, 2002.

\bibitem[Jin-Xin and Zhou(2017)]{Jin2017On}
Jin-Xin and Zhou.
\newblock On g-extra connectivity of hypercube-like networks.
\newblock \emph{Journal of Computer \& System Sciences}, 2017.

\bibitem[Li and Yang(2013)]{2013Bounding}
H.~Li and W.~Yang.
\newblock Bounding the size of the subgraph induced by m vertices and extra
  edge-connectivity of hypercubes.
\newblock \emph{Discrete Applied Mathematics}, 161\penalty0 (16-17):\penalty0
  2753--2757, 2013.

\bibitem[Sampathkumar(1984)]{1984Connectivity}
E.~Sampathkumar.
\newblock Connectivity of a graph - a generalization.
\newblock \emph{J.combin.inform.system Sci}, 1984.

\bibitem[Shang et~al.(2019)Shang, Sabir, Meng, and Guo]{2019Characterizations}
H.~Shang, E.~Sabir, J.~Meng, and L.~Guo.
\newblock Characterizations of optimal component cuts of locally twisted cubes.
\newblock \emph{The Bulletin of the Malaysian Mathematical Society Series 2},
  43\penalty0 (3):\penalty0 1--17, 2019.

\bibitem[Shen et~al.(2015)Shen, Qiu, and Cheng]{2015Connectivity}
Z.~Shen, K.~E. Qiu, and E.~Cheng.
\newblock Connectivity results of complete cubic networks as associated with
  linearly many faults.
\newblock \emph{Journal of interconnection networks}, 2015.

\bibitem[Zhao et~al.(2018)Zhao, Shuli, Yang, Weihua, Zhang, Shurong, Xu, and
  Liqiong]{Zhao2018Component}
Zhao, Shuli, Yang, Weihua, Zhang, Shurong, Xu, and Liqiong.
\newblock Component edge connectivity of hypercubes.
\newblock \emph{International Journal of Foundations of Computer Science},
  2018.

\bibitem[Zhao and Yang(2019)]{2019Conditional}
S.~Zhao and W.~Yang.
\newblock Conditional connectivity of folded hypercubes.
\newblock \emph{Discrete Applied Mathematics}, 257:\penalty0 388--392, 2019.

\end{thebibliography}

%\end{CJK*}

\end{document}